\documentclass[11pt]{article}
\usepackage{latexsym,amsmath,stackrel,color,amsthm,amssymb,epsfig,graphicx,mathrsfs}
\usepackage{graphicx}
\usepackage{amssymb}
\usepackage[left=1in,top=1in,right=1in,bottom=1in]{geometry}
\usepackage[linktocpage=true]{hyperref}
\usepackage{setspace}
\usepackage{amssymb, amsmath, amsthm, graphicx,mathrsfs}
\usepackage{caption}
\usepackage{titlesec}

\usepackage{comment,caption}
\usepackage{tikz}
\makeatletter
\def\thm@space@setup{%
  \thm@preskip=\parskip \thm@postskip=0pt
}
\makeatother

\usepackage{thmtools}

\declaretheoremstyle[%
  spaceabove=6pt,%
  spacebelow=6pt,%
  headfont=\normalfont\itshape,%
  postheadspace=1em,%
  qed=\qedsymbol%
]{mystyle}

\def\qed{\hfill\ifhmode\unskip\nobreak\fi\quad\ifmmode\Box\else\hfill$\Box$\fi}
\def\ite#1{\hfill\break${}$\hbox to 50pt {\quad(#1)\hfill}}
\newtheorem{thm}{Theorem}
\newtheorem{cor}[thm]{Corollary}

\newtheorem{lemma}[thm]{Lemma}

\tikzstyle{vertex}=[circle,fill=black,inner sep=2pt]
\tikzstyle{vertrect}=[draw,rectangle,inner sep=2pt]
\tikzstyle{vertdia}=[draw,diamond,inner sep=2pt]
\newcommand{\vb}[1]{\boldsymbol{#1}}

\newcommand{\Z}{\mathbb{Z}}

\newcommand{\A}{\mathcal{A}}
\newcommand{\M}{\mathcal{M}}
\newcommand{\h}{\mathcal{H}}

\titleformat{\subsection}[runin]% runin puts it in the same paragraph
        {\normalfont\bfseries}% formatting commands to apply to the whole heading
        {\thesubsection}% the label and number
        {0.4em}% space between label/number and subsection title
        {}% formatting commands applied just to subsection title
        [.]% punctuation or other commands following subsection title

\begin{document}
\title{A generalization of the Bollob{\'a}s set pairs inequality}
\author{
	Jason O'Neill \\
	Department of Mathematics, \\
	University of California, San Diego \\
	\texttt{jmoneill@ucsd.edu}
	\and
	Jacques Verstra\"{e}te
\footnote{Research supported by NSF award DMS-1800332} \\
	Department of Mathematics, \\
	University of California, San Diego \\
	\texttt{jverstra@math.ucsd.edu}
}
\maketitle

\vspace{-0.1in}

\begin{abstract}
The Bollob{\'a}s set pairs inequality is a fundamental result in extremal set theory with many applications. In this paper, for $n \geq k \geq t \geq 2$, we consider a collection of $k$ families $\mathcal{A}_i: 1 \leq i \leq k$ where $\mathcal{A}_i = \{ A_{i,j} \subset [n] : j \in [n] \}$ so that $A_{1, i_1} \cap \cdots \cap A_{k,i_k} \neq  \emptyset $ if and only if there are at least $t$ distinct indices $i_1,i_2,\dots,i_k$. Via a natural connection to a hypergraph covering problem, we give bounds on the maximum size $\beta_{k,t}(n)$ of the families with ground set $[n]$.
\end{abstract}

\vspace{-0.1in}

\section{Introduction}

A central topic of study in extremal set theory is the maximum size of a family of subsets of an $n$-element set subject to restrictions on their intersections. Classical theorems in the area are discussed in Bollob\'{a}s~\cite{B}. In this paper, we generalize one such theorem, known as the Bollob\'{a}s set pairs inequality or two families theorem~\cite{BBOL}:

\begin{thm} {\rm (Bollob{\'a}s)} \label{thm:Bollobas}
Let $\mathcal{A} = \{A_1,A_2,\dots,A_m\}$ and $\mathcal{B} = \{B_1,B_2,\dots,B_m\}$ be families of finite sets, such that $A_i \cap B_j \neq \emptyset$ if and only if $i,j \in [m]$ are distinct.  Then
\begin{equation}\label{setpairs}
\sum_{i = 1}^m  \binom{|A_i \cup B_i|}{|A_i|}^{-1} \leq 1.
\end{equation}
\end{thm}

For convenience, we refer to a pair of families $\mathcal{A}$ and $\mathcal{B}$ satisfying the conditions of Theorem \ref{thm:Bollobas} as a {\em Bollob\'{a}s set pair}.  The inequality above is tight, as we may take the pairs $(A_i,B_i)$ to be distinct partitions of a set of size $a + b$ with $|A_i| = a$ and $|B_i| = b$ for $1 \leq i \leq {a + b \choose b}$. The latter inequality was proved for $a = 2$  by Erd\H{o}s, Hajnal and Moon~\cite{EHM}, and in general has a number of different proofs~\cite{H,JP,Katona,L1,L2}. A geometric version was proved by Lov\'{a}sz~\cite{L1,L2}, who showed that if $A_1,A_2,\dots,A_m$ and $B_1,B_2,\dots,B_m$ are respectively $a$-dimensional and $b$-dimensional subspaces of a linear space and $\mbox{dim}(A_i \cap B_j) = 0$ if and only if $i,j \in [m]$ are distinct, then $m \leq {a + b \choose a}$.

\subsection{Main Theorem}
Theorem \ref{thm:Bollobas} has been generalized in a number of different directions in the literature~\cite{Fr,F,K,L,JTOL,T1}.
In this paper, we give a generalization of Theorem \ref{thm:Bollobas} from the case of two families to $k \geq 3$ families of sets with conditions on the $k$-wise intersections. For $2 \leq t \leq k$, a {\em Bollob\'{a}s $(k,t)$-tuple} is a sequence $(\mathcal{A}_1,\mathcal{A}_2,\dots,\mathcal{A}_k)$ of set families  $\mathcal{A}_{j} = \{A_{j,i} : 1 \leq i \leq m\}$ where
$\bigcap_{j = 1}^k A_{j,i_j} \neq \emptyset$ if and only if at least $t$  of the indices $i_1,i_2,\dots,i_k$ are distinct. We refer to $m$ as the {\em size} of the
Bollob\'{a}s $(k,t)$-tuple. Let $[m]_{(t)}$ denote the set of sequences of $t$ distinct elements of $[m]$ and fix a surjection $\phi : [k] \rightarrow [t]$.
For $\sigma \in [m]_{(t - 1)}$, set $\sigma(t) = \sigma(1)$ and define $A_{1,\sigma}(\phi) = \bigcap_{j : \phi(j) = 1} A_{j,\sigma(1)}$ and, for $2 \leq j \leq t$, we define
\[ A_{j,\sigma}(\phi) = \bigcap_{h : \phi(h) = j} A_{h,\sigma(j)} \backslash \bigcup_{h = 1}^{j - 1} A_{h,\sigma}(\phi).\]
Using this notation, we generalize (\ref{setpairs}) as follows:

\begin{thm}\label{thm:genthreshk}
Let $k \geq t \geq 2$ and $m \geq t$, let $\phi : [k] \rightarrow [t]$ be a surjection, and let
$(\mathcal{A}_1,\mathcal{A}_2,\dots,\mathcal{A}_k)$ be a Bollob{\'a}s $(k,t)$-tuple of size $m$. Then
\begin{equation}\label{eq:bskicasekequalsk}
\sum_{\sigma \in [m]_{(t - 1)}} {|A_{1,\sigma}(\phi) \cup A_{2,\sigma}(\phi) \cup \cdots \cup A_{t,\sigma}(\phi)| \choose |A_{1,\sigma}(\phi)| \;\; |A_{2,\sigma}(\phi)| \;\; \cdots \;\; |A_{t,\sigma}(\phi)|}^{-1} \leq 1.
\end{equation}
\end{thm}

We show in Section \ref{subsec:sharpness} that this inequality is tight for all $k \geq t = 2$, but do not have an example to show that this inequality is tight for any $t > 2$. 

For $n \geq k \geq t \geq 2$, let $\beta_{k,t}(n)$ denote the maximum $m$ such that there exists a Bollob\'{a}s $(k,t)$-tuple
of size $m$ consisting of subsets of $[n]$. Then (\ref{setpairs}) gives $\beta_{2,2}(n) \leq {n \choose \lfloor n/2 \rfloor}$ which is tight for all $n \geq 2$. Letting  $H(q) = -q \log_2 q - (1-q) \log_2(1-q)$ denote the standard binary entropy function, we prove the following theorem:

\begin{thm}\label{thm:betabounds}
For $k \geq 3$ and large enough $n$,
\begin{equation}\label{eq:betaboundstequals2}
\frac{1} {k} \leq \frac{\log_2 \beta_{k,2}(n)}{n} \leq H\left(\frac{1}{k}\right) \leq \frac{\log_2 (ke)}{k}.
\end{equation}
For $k \geq t \geq 3$ and large enough $n$,
\begin{equation}\label{eq:betaboundsgent}
\frac{\log_2 e}{\binom{k}{t-1}(t+1)t^{t-1}} \leq  \frac{\log_2 \beta_{k,t}(n)}{n} \leq \frac{2}{\binom{k}{t-1}(t-1)^{t-3}}.
\end{equation}
\end{thm}

This determines $\log_2 \beta_{k,2}(n)$ up to a factor of order $\log_2 k$ and $\log_2 \beta_{k,2}(n)$ up to a factor of order $t$. We leave it as an open problem to determine the asymptotic value of
$(\log_2 \beta_{k,t}(n))/n$ as $n \rightarrow \infty$ for any $k \geq 3$ and $t \geq 2$. A natural source for lower bounds on $\beta_{k,t}(n)$ comes
from the probabilistic method -- see the random constructions in Section \ref{subsec:ubfkt} which establish the lower
bounds in Theorem \ref{thm:betabounds}. To prove Theorem \ref{thm:betabounds}, we use a natural connection to hypergraph covering problems.

\subsection{Covering hypergraphs}
Theorem \ref{thm:Bollobas} has a wide variety of applications, from saturation problems~\cite{BBOL,MS} to covering problems for graphs~\cite{H,Orlin}, complexity of 0-1 matrices~\cite{Tarjan}, geometric problems~\cite{AK}, counting cross-intersecting families~\cite{FrK}, crosscuts and transversals of hypergraphs~\cite{T1,T2,T3},
hypergraph entropy~\cite{KM,Simonyi}, and perfect hashing~\cite{FK,GR}. In this section, we give an application of our main results to hypergraph covering problems. For a $k$-uniform hypergraph $H$, let $f(H)$ denote the minimum number of complete $k$-partite $k$-uniform hypergraphs whose union is $H$. In the case of graph covering, a simple connection to the Bollob\'{a}s set pairs inequality (\ref{setpairs})
may be described as follows. Let $H$ denote the complement of  a perfect matching $\{x_i y_i : 1 \leq i \leq n\}$ in the complete bipartite graph $K_{n,n}$ with parts $X = \{x_1,x_2,\dots,x_n\}$ and $Y = \{y_1,y_2,\dots,y_n\}$. If $H_1,H_2,\dots,H_m$ are complete bipartite graphs in a minimum covering of $H$, then let $A_i = \{j : x_i \in V(H_j)\}$ and $B_i = \{j : y_i \in V(H_j)\}$. Setting $\mathcal{A}=\{A_i\}_{i \in [m]}$ and $\mathcal{B}=\{B_i\}_{i \in [m]}$, it is straightforward to check that $(\mathcal{A}, \mathcal{B})$ is a Bollob\'{a}s set pair, and Theorem \ref{thm:Bollobas} applies to give
\begin{equation}\label{bicliquecover}
f(K_{n,n} \backslash M) = \min\{m : {m \choose \lceil m/2\rceil} \geq n\}.
\end{equation}
In a similar way, Theorem \ref{thm:genthreshk} applies to covering complete $k$-partite $k$-uniform hypergraphs. Let $K_{n,n,\dots,n}$ denote the complete $k$-partite $k$-uniform hypergraph with parts $X_{i} = \{x_{ij} : j \in [n]\}$
for $i \in [k]$. Let $H_{k,t}(n)$ denote the subhypergraph consisting of hyperedges $\{x_{1,i_1},x_{2,i_2},\dots,x_{k,i_k}\}$ such that at least $t$ of the indices $i_1,i_2,\dots,i_k$ are distinct, and set $f_{k,t}(n) = f(H_{k,t}(n))$. Then there is a one-to-one correspondence between Bollob\'{a}s $(k,t)$-tuples of subsets of $[m]$ and coverings of $H_{k,t}(n)$ with $m$ complete $k$-partite $k$-graphs. We let $\beta_{k,t}(m)$ be the maximum size of a Bollob\'{a}s $(k,t)$-tuple of subsets of $[m]$, so that
\begin{equation}\label{connection}
f_{k,t}(n) = \min\{m : \beta_{k,t}(m) \geq n\}.
\end{equation}
This correspondence together with Theorem \ref{thm:genthreshk} will be exploited to prove
\begin{equation}\label{nicebound}
f_{k,2}(n) \geq \min\{m : {m \choose \lceil m/k \rceil} \geq n \}
\end{equation}
which is partly an analog of (\ref{bicliquecover}). More generally, we prove the following theorem:

\begin{thm}\label{thm:generalktbounds} 
For $k \geq 3$ and large enough $n$,
\begin{equation}\label{eq:coveringboundstequals2}
\frac{k}{\log_2 (ke)} \leq \frac{1}{H(\frac{1}{k})} \leq \frac{f_{k,2}(n)}{\log_2 n} \leq k.
\end{equation}
For $k \geq t \geq 3$ and large enough $n$,
\begin{equation}\label{eq:generalbounds}
\binom{k}{t-1} \frac{(t-1)^{t-3}}{2} \leq \frac{f_{k,t}(n)}{\log_2n} \leq \frac{(t+1)t^{t-1}}{\log_2 e} \binom{k}{t-1}.
\end{equation}
\end{thm}

The bounds on $\beta_{k,t}(n)$ in Theorem (\ref{thm:betabounds}) follow immediately from this theorem and (\ref{connection}). Equation \eqref{eq:generalbounds} gives the order of magnitude for each $t \geq 3$ as $k \to \infty$, but for $t=2$, Equation \eqref{eq:coveringboundstequals2} has a gap of order $\log_2 k$. From (\ref{nicebound}), we obtain $\beta_{k,2}(n) \leq {n \choose \lfloor n/k\rfloor}$.
It is perhaps unsurprising that the asymptotic value of $f_{k,t}(n)/\log_2 n$ as $n \rightarrow \infty$ is not known for any $k > 2$, since a limiting value of $f(K_n^k)/\log_2 n$ is not known for any $k > 2$ -- see K\"{o}rner and Marston~\cite{KM} and Guruswami and Riazanov~\cite{GR}.

\subsection{Organization and notation}
This paper is organized as follows. In Section \ref{mainproof}, we prove Theorem \ref{thm:genthreshk}. In Section \ref{subsec:sharpness}, we construct a Bollob{\'a}s $(k,2)$-tuple which achieves equality in Theorem \ref{thm:genthreshk} and in Section \ref{subsec:eck3t2}, we construct a Bollob{\'a}s $(k,2)$-tuple which gives the lower bound in Equation \eqref{eq:betaboundstequals2}. The upper bound on $f_{k,t}(n)$ in Theorem \ref{thm:generalktbounds} comes from a probabilistic construction in Section \ref{subsec:ubfkt}, and the proof of the lower bound on $f_{k,t}(n)$ is given in Section \ref{subsection:lbfkk}; we prove (\ref{nicebound}) 
in Section \ref{subsec:lbfk2}.

\section{Proof of Theorem \ref{thm:genthreshk}}\label{mainproof}
Given a Bollob{\'a}s set $(k,t)$-tuple $(\A_1, \ldots, \A_k)$ with $\A_j= \{A_{j,i} : 1 \leq i \leq m\}$ and a surjection $\phi:[k] \to [t]$, consider $\A_{\ell}(\phi):1\leq \ell \leq t$ where $\A_\ell(\phi) = \{A_{\ell,i}(\phi) : 1 \leq i \leq m\} $ and
\[ A_{\ell,i}(\phi) = \bigcap_{h : \phi(h) = \ell} A_{h,i}.\]
It follows that $(\A_1(\phi), \ldots, \A_t(\phi))$ is a Bollob{\'a}s set $(t,t)$-tuple and hence it suffices to prove Theorem \ref{thm:genthreshk} in the case where $t=k$. In this setting, surjections $\phi:[k] \to [k]$ simply permute the $k$ families and as such we suppress the notation of $\phi$ for the remainder of this section.    

Let $(\A_1, \ldots, \A_k)$ with $\A_j= \{A_{j,i} : 1 \leq i \leq m\}$ be a Bollob{\'a}s set $(k,k)$-tuple, and ground set 
\[ X = \bigcup_{i=1}^m (A_{1,i} \cup A_{2,i}\cup \cdots \cup A_{k,i})\]
with $|X|=n$. For $\sigma \in [m]_{(k-1)}$, define a subset $\mathscr{C}_{\sigma}$ of permutations $\pi:X \to [n]$ by \[ \mathscr{C}_{\sigma}:=
\left\{\pi: X \to [n]: \max_{ x \in A_{1,\sigma}} \pi(x) < \min_{ y \in A_{2,\sigma} }  \pi(y) \leq \max_{ y \in A_{2,\sigma}} \pi(y) < \cdots < \min_{ z \in A_{k,\sigma}} \pi(z) \right\}. \]

Letting $U_\sigma:= A_{1,\sigma} \cup \cdots \cup A_{k,\sigma}$, elementary counting methods give
\begin{equation}\label{eq:countingchains}
|\mathscr{C}_{\sigma}| = \binom{n}{|U_\sigma|} \, |A_{1,\sigma}|! \cdots |A_{k,\sigma}|! (n- |U_\sigma|)! = n! \cdot \binom{|U_\sigma|}{|A_{1,\sigma}| \cdots |A_{k,\sigma}|}^{-1}.    
\end{equation}

We will now prove a lemma which states that $\{\mathscr{C}_\sigma \}_{\sigma \in [m]_{(k-1)}}$ forms a disjoint collection of a permutations. The general proof only works for $k \geq 4$, so we first consider $k=3$. 

\begin{lemma}\label{lem:5cases}
If $\sigma_1, \sigma_2 \in [m]_{(2)}$ are distinct, then $\mathscr{C}_{\sigma_1} \cap \mathscr{C}_{\sigma_2} = \emptyset$.	
\end{lemma}

\begin{proof}
Seeking a contradiction, suppose there exists $\pi \in \mathscr{C}_{\sigma_1} \cap \mathscr{C}_{\sigma_2}.$ After relabeling, it suffices to consider the following five cases. 
\begin{center}
\begin{tabular}{ll}
(1) $\sigma_1= \{1,3\}$ and $\sigma_2 = \{2,4\}$ & (2) $\sigma_1= \{1,3\}$ and $\sigma_2 = \{2,3\}$ \\ (3) $\sigma_1= \{1,2\}$ and $\sigma_2 = \{1,3\}$ &
(4)	$\sigma_1= \{1,2\}$ and $\sigma_2 = \{2,3\}$ \\  (5) $\sigma_1= \{1,2\}$ and $\sigma_2 = \{3,1\}$.
\end{tabular}
\end{center}

In case (1), without loss of generality, $\max\{ \pi(x) : x  \in A_{1,1}\} \leq \max\{ \pi(x): x  \in A_{1,2}\}$ and thus $\pi \in \mathscr{C}_{\sigma_2}$ yields  \[ \ \max_{x \in A_{1,1}} \pi(x) \leq \max_{x \in A_{1,2}} \pi(x) < \min_{ y \in A_{2,4} \setminus A_{1,2}} \pi (y).\] Then as $A_{1,1} \cap A_{2,4} \cap A_{3,2} \neq \emptyset$, there exists $w \in A_{1,1} \cap A_{2,4} \cap A_{3,2}$. It follows that $w \notin A_{1,2}$ since if $w \in A_{1,2}$, then $w \in A_{1,2} \cap A_{2,4} \cap A_{3,2} \neq \emptyset$; a contradiction. But this yields a contradiction as \[ \pi(w) \leq \max_{x \in A_{1,1}} \pi(x) \leq \max_{x \in A_{1,2}} \pi(x) < \min_{ y \in A_{2,4} \setminus A_{1,2}} \pi (y) \leq \pi(w). \] 

In case (2), without loss of generality, $\max \{ \pi(x) : x \in A_{1,1}\} \leq \max\{ \pi(x) : x \in A_{1,2}\} $ and we recover a similar contradiction as case (1) by noting that there exists $w \in A_{1,1} \cap A_{2,3} \cap A_{3,2}$ with $w \notin A_{1,2}$.

In case (3), without loss of generality, $\max\{ \pi(x) : x \in A_{2,2} \setminus A_{1,1}\}  \leq \max\{ \pi(x) : x \in A_{2,3} \setminus A_{1,1} \} $ and $\pi \in \mathscr{C}_{1,3}$ yields $\max \{ \pi(x) :  x \in A_{2,3}  \setminus A_{1,1} \} < \min\{ \pi(x) : x \in A_{3,1} \setminus (A_{1,1} \cup A_{2,3})\}.$ Thus \[ \max\{ \pi(x) : x \in A_{2,2} \setminus A_{1,1}\} < \min\{ \pi(x) : x \in A_{3,1} \setminus (A_{1,1} \cup A_{2,3})\} \] and there exists $w \in A_{1,3} \cap A_{2,2} \cap A_{3,1}$ with $ w \notin A_{1,1}$ and $w \notin A_{2,3}$. It follows that $\pi(w) < \pi(w)$, a contradiction.

In case (4), if $\max\{\pi(x) : x \in A_{1,1}\} \leq \max\{ \pi(x) : x \in A_{1,2}\}$, then using $w \in A_{1,1} \cap A_{2,3} \cap A_{3,2}$ and noting $w \notin A_{1,2}$, we get a contradiction. Thus, we may assume otherwise and $\pi \in \mathcal{C}_{1,2}$ gives \[ \max_{x \in A_{1,2}} \pi(x) < \max_{x \in A_{1,1}} \pi(x) < \min_{z \in A_{3,1} \setminus (A_{1,1} \cup A_{2,2})} \pi(z). \] This is a contradiction as there exists $w \in A_{1,2} \cap A_{2,3} \cap A_{3,1}$ with $ w \notin A_{1,1}$ and $w \notin A_{2,2}$.

In case (5), if $\max\{ \pi(x): x \in A_{1,1}\} \leq \max\{ \pi(x) : x \in A_{1,3}\}$, then we may proceed as in the latter part of case (4) using $w \in A_{1,1} \cap A_{2,2} \cap A_{3,3}$ and $w \notin A_{2,1}$ and $w \notin A_{1,3}$ to get a contradiction. Otherwise, proceeding as in case (1) and noting there exists $w \in A_{1,3} \cap A_{2,2} \cap A_{3,1}$, but $w \notin A_{1,1}$ yields a contradiction. \qedhere  

\end{proof}

A similar argument yields the analog of Lemma \ref{lem:5cases} to the case where $k \geq 4$.

\begin{lemma} \label{lemma:distinctchainfam}
Let $k \geq 4$. If $\sigma_1, \sigma_2 \in [m]_{(k-1)}$ are distinct, then $\mathscr{C}_{\sigma_1} \cap \mathscr{C}_{\sigma_2} = \emptyset$.
\end{lemma}

\begin{proof}
Since $\sigma_1, \sigma_2 \in [m]_{(k-1)}$ are distinct, there exists minimal $h \in [k-1]$ so that $\sigma_1(h) \neq \sigma_2(h)$. Seeking a contradiction, suppose there exists a $\pi \in \mathscr{C}_{\sigma_1} \cap \mathscr{C}_{\sigma_2}$.  Without loss of generality, \[ \max\{ \pi(x) : x \in A_{h, \sigma_1}\} \leq \max\{ \pi(x) : x \in A_{h,\sigma_2}\}  < \min\{ \pi(z) : z \in A_{k, \sigma_2}\}. \] Now, consider a bijection $\tau: [k-1] \setminus \{h\} \to [k-1] \setminus \{1\}$ which has no fixed points. As in Lemma \ref{lem:5cases}, we want to show that there exists a $w \in A_{h,\sigma_1} \cap A_{k, \sigma_2}$ and consider two separate cases.

\medskip
First, suppose that $\sigma_1(h) \notin \sigma_2 ([k-1])$. As $|\{\sigma_1(h), \sigma_2(1), \ldots, \sigma_2(k-1)\}| =k$, there exists
\begin{equation}\label{eq: contrafirstcase}
w \in A_{h,\sigma_1(h)} \cap A_{k, \sigma_2(1)} \cap \bigcap\limits_{l \in [k-1] \setminus \{h\}} A_{l,\sigma_2(\tau(l))}.
\end{equation}

Next, suppose that $\sigma_1(h) = \sigma_2(x)$ for some $x$. We now claim that $x\neq 1$. If $h=1$, then this is trivial. If $h>1$, then $\sigma_1(1)= \sigma_2(1)$, so $\sigma_1(h) \neq \sigma_2(1)$ since $\sigma_1(h) \neq \sigma_1(1)$. For $\tau$ as above, there exists $y \in [k-1] \setminus \{h\}$ so that $ \tau(y) = x $. Taking $\gamma$ distinct from $\{ \sigma_2(1), \ldots, \sigma_2(k-1)\} \setminus \{\sigma_2(x)\}$, $|\{ \sigma_1(h), \gamma, \sigma_2(1), \ldots, \sigma_2(k-1) \} \setminus \{\sigma_2(x)\} | = k $ and hence there exists
\begin{equation}\label{eq: contrasecondcase}
w \in A_{h,\sigma_1(h)} \cap A_{k,
\sigma_2(1)} \cap A_{y,\gamma} \cap \bigcap\limits_{l \in [k-1] \setminus \{y,h\}} A_{l,\sigma_2(\tau(l))}.
\end{equation}

By construction, $w \in A_{h, \sigma_1(h)} \cap A_{k,\sigma_2(1)}$. Suppose there exists a $t \in [k-1] \setminus \{h\}$ so that $ w\in A_{t,\sigma_2(t)}$. As $\tau$ has no fixed points, replacing the set in the $k$-wise intersection corresponding to $\A_t$ with $A_{t,\sigma_2(t)}$  in either \eqref{eq: contrafirstcase} or \eqref{eq: contrasecondcase}, $w$ is an element of this new $k$-wise intersection with $(k-1)$ distinct indices; a contradiction. If $w \in A_{h, \sigma_2(h)}$, then we may similarly replace $A_{h, \sigma_1(h)}$ with $A_{h ,\sigma_2(h)}$ in the $k$-wise intersection in either \eqref{eq: contrafirstcase} or \eqref{eq: contrasecondcase} to get a contradiction. Thus, $w \notin  A_{1,\sigma_2(1)} \cup \cdots \cup A_{k-1,\sigma_2(k-1)}$ and hence $w \in A_{h,\sigma_1} \cap A_{k, \sigma_2}$ so that $\pi(w)<\pi(w)$; a contradiction. 
\end{proof}

Using Equation \eqref{eq:countingchains}, Lemma \ref{lem:5cases}, and Lemma \ref{lemma:distinctchainfam}, we are now able to prove Theorem \ref{thm:genthreshk} in the case where $t=k$.

There are $n!$ total permutations, and Lemma \ref{lem:5cases} and Lemma \ref{lemma:distinctchainfam} yield that each of which appears in at most one of the sets $\mathscr{C}_{\sigma}$ for $\sigma \in [m]_{(k-1)}$. Hence, using $|\mathscr{C}_{\sigma}|$ in Equation \eqref{eq:countingchains}, \[\sum_{ \sigma \in [m]_{(k-1)} } |\mathscr{C}_\sigma| = \sum_{\sigma \in [m]_{(k-1)}} n! \cdot \binom{|A_{1,\sigma}  \cup \cdots \cup A_{k,\sigma}| }{|A_{1, \sigma}| \cdots  |A_{k, \sigma}|}^{-1} \leq n! \] and thus the result follows by dividing through by $n!$.

\subsection{Sharpness of Theorem \ref{thm:genthreshk}}\label{subsec:sharpness}
We give a simple construction establishing the sharpness of Theorem \ref{thm:genthreshk} for $k \geq t = 2$. Let $n \geq 4k$ and using addition modulo $n$, define $A_{1,i} = \{i\}^c$, $A_{j,i} = \{ i-(j-1), i+(j-1)\}^c$ for $j \in [2,k-1]$, and $A_{k,i} = \{ i-k+2, i-k+3, \ldots, i+k-2 \}$. Letting $\mathcal{A}_j =\{A_{j,i}\}_{i\in [n]}$ for all $j \in [k]$, we will show $(\mathcal{A}_1, \ldots, \mathcal{A}_k)$ is a Bollob{\'a}s  $(k,2)$-tuple. 
Since $|A_{1,i}| = n-1$ and $|A_{2,i} \cap \cdots \cap A_{k,i}|=1$, Theorem \ref{thm:genthreshk} with $t = 2$ gives
\[  1 \;\;\geq \;\; \sum_{i=1}^{n} \binom{|A_{1,i}| + |A_{2,i} \cap \cdots \cap A_{k,i}| }{|A_{1,i}|} ^{-1} \; = \; \; \sum_{i=1}^{n} \frac{1}{n} \; \;= \; \; 1. \]

By construction, for all $i \in [n]$, $A_{1,i}\cap A_{2,i} \cap \cdots \cap A_{k,i} = \emptyset$. It thus suffices to show these are the only empty $k$-wise intersections. To this end, for $\vb{i} = (i_1, \ldots, i_{k-1})$, define \[ A(\vb{i}):= (A_{1,i_1} \cap \cdots \cap A_{k-1,i_{k-1}})^c. \] 

\begin{lemma}\label{lemma: examplethres2}
Let $\vb{i} = (i_1, \ldots, i_{k-1})$. If $A(\vb{i})^c= A_{k, i_k}$, then  $i_1 = \cdots = i_k$.
\end{lemma}

\begin{proof}
We proceed by induction on $k$ where the result is trivial when $k=2$. In the case where $k>2$, $i_{k-1}-k+2 = i_k +x$ for some $x$ such that $-(k-2) \leq x \leq (k-2)$ and thus $i_{k-1} +(k-2) = i_{k-1} - (k-2) + (2k-4) = i_k+x +(2k-4).$

Next, there is a $y$ such that  $-(k-2) \leq y \leq (k-2)$ with $ i_{k-1} +(k-2) = i_k + y$, and since $n \geq 4k$, $x+2k-4 = y$ with equality over $\Z$ and moreover $i_{k-1}+(k-2)= i_k+(k-2)$ over $\Z$ and hence $i_k=i_{k-1}$. Removing these elements from each set, the result then follows by induction.
\end{proof}

If $A_{1,i_1} \cap \cdots \cap A_{k,i_k} = \emptyset$, then as
$A(\vb{i}) = A_{1,i_1} \cap A_{2,i_2} \cap \cdots A_{k-1,i_{k-1}}$,
\[ \emptyset = A_{1,i_1} \cap A_{2,i_2} \cap \cdots \cap  A_{k-1,i_{k-1}} \cap A_{k,i_k} = A(\vb{i}) \cap A_{k,i_k}.\]

The result follows by examining the cardinality of $A(\vb{i})$ and $A_{k,i_k}$ and using Lemma \ref{lemma: examplethres2}. 	

\subsection{An Explicit Construction}\label{subsec:eck3t2} 
The constructions of large Bollob\'{a}s $(k,t)$-tuples of subsets of $[n]$ are probabilistic -- see Section \ref{subsec:ubfkt}.
An explicit construction of a Bollob{\'a}s $(3,2)$-tuple $(\mathcal{A}_1,\mathcal{A}_2,\mathcal{A}_3)$ where $|\mathcal{A}_i| = 2^n$ and each $\mathcal{A}_i$ consists of subsets of $X$ for $|X|=3n$ may be described as follows. Let $I_j:=\{ x_{j,1}, x_{j,2}, x_{j,3} \}$ and consider $X = I_1 \sqcup \cdots \sqcup I_n$. Now, for each $f:[n] \to [3]$, define 
\begin{eqnarray*}
A_{1,f} &:=& \{ x_{1,f(1)}, \; \ldots \;, x_{n,f(n)} \}^c\\
A_{2,f} &:=& \{ x_{1,f(1)+1}, \;  \ldots \; , x_{n,f(n)+1} \}^c \\
A_{3,f} &:=& \{ x_{1,f(1)+2}, \; \ldots \; , x_{n,f(n)+2} \}^c
\end{eqnarray*}
where we work modulo $3$ within the subscripts of $I_j$. Letting $\A_i:=\{ A_{i,f} : f:[n] \to [2] \}$, for $i \in [3]$,  it is straightforward to check that $(\mathcal{A}_1,\mathcal{A}_2,\mathcal{A}_3)$ is a Bollob{\'a}s $(3,2)$-tuple.

The above construction generalizes for $k>3$. An explicit construction of a Bollob{\'a}s $(k,2)$-tuple $(\mathcal{A}_1,\mathcal{A}_2, \ldots, \mathcal{A}_k)$ where $|\mathcal{A}_i| = 2^n$ and each $\mathcal{A}_i$ consists of subsets of $X$ for $|X|=kn$ may be described as follows. Let $I_j:=\{ x_{j,1}, x_{j,2}, \ldots,  x_{j,k} \}$ and consider $X = I_1 \sqcup \cdots \sqcup I_n$. Now, for each $f:[n] \to [2]$ and $j \in [k]$, define 
\[ A_{j,f} := \{ x_{1,f(1)+j-1}, \; \ldots \;, x_{n,f(n)+j-1} \}^c \] 
where we work modulo $k$ within the subscripts of $I_j$. It is straightforward to check that $(\mathcal{A}_1,\mathcal{A}_2,\ldots, \mathcal{A}_k)$ is a Bollob{\'a}s $(k,2)$-tuple. This establishes the lower bound on $\beta_{k,2}(n)$ in Equation \eqref{eq:betaboundstequals2} and hence the upper bound on $f_{k,2}(n)$ in Equation \eqref{eq:coveringboundstequals2}.
 
\section{Proof of Theorem \ref{thm:generalktbounds}}\label{sec:generalktbounds}
\subsection{Upper bound on $f_{k,t}(n)$}\label{subsec:ubfkt}

We wish to find a covering of $H_{k,t}(n)$ with complete $k$-partite $k$-graphs and assume the parts of $H_{k,t}(n)$ are
$X_1,X_2,\dots,X_k$. For each subset $T$ of $[k]$ of size $t$, consider the uniformly randomly coloring $\chi_T: [n] \to T$. Given such a $\chi_T$, let $Y_i \subset X_i$ be the vertices of color $i$ for $i \in T$; that is $Y_i := \{ x_{ij} : \chi(j) = i   \}$ and $Y_i = X_i$ for $i \notin T$. Denote by $H(T,\chi)$ the (random) complete $k$-partite hypergraph with parts $Y_1,Y_2,\dots,Y_k$, and note that $H(T, \chi) \subset H_{k,t}(n)$. We place each $H(T,\chi)$ a total of $N$ times independently and randomly where 
\[ N = \Big\lfloor \frac{(t+1)t^t \log_2 n}{(k-t+1) \log_2 e} \Big\rfloor \]
and produce ${k \choose t}N$ random subgraphs $H(T,\chi)$. For a set partition $\pi$ of $[k]$, let $|\pi|$ denote the number of parts in the partition and index the parts by $[|\pi|]$. Given a set partition  $\pi = (P_1,P_2,\dots,P_s)$, let
\[ f(\pi, t) = \sum_{T \in [s]^{(t)}} \prod_{i \in T} |P_i|. \]
If $U$ is the number of edges of $H_{k,t}(n)$ not in any of these subgraphs, then

\begin{equation}\label{eq:expectationkt}
\mathbb E(U) \leq \sum_{|\pi| \geq t} n^{|\pi|} (1 - t^{-t})^{N f(\pi,t)} = \sum_{t \leq s \leq k} n^s \sum_{|\pi|=s} (1 - t^{-t})^{N f(\pi,t)}.
\end{equation}

For sufficiently large $n$, we claim that $\mathbb E(U) < 1$, which implies there exists a covering of $H_{k,t}(n)$ with at most ${k \choose t}N$
complete $k$-partite $k$-graphs, as required. The following technical lemma states that $f$ is a decreasing function in the set partition lattice, and
that $f(\pi,t)$ increases when we merge all but one element of a smaller part of $\pi$ with a larger part of $\pi$:

\begin{lemma}\label{lemma:setpartitionstat}
Let $k \geq s \geq t \geq 2$, and let $\pi = (P_1,P_2,\dots,P_s)$ be a partition of $[k]$.
\begin{center}
\begin{tabular}{lp{5.8in}}
{\rm (i)} & If $\pi'$ is a refinement of $\pi$ with $|\pi'|=s+1$, then  $f(\pi,t) \leq f(\pi',t)$. \\
{\rm (ii)} & If $|P_1| \geq |P_2| \geq 2$ and $a \in P_2$, and $\pi'$ is the partition $(P_1',P_2',\dots,P_s')$ of $[k]$ with $P_1' = P_1 \cup P_2 \setminus \{a\}$ and $P_2' = \{a\}$ and with $P_i' = P_i$ for $3 \leq i \leq s$, then $f(\pi',t) \leq f(\pi, t)$.	
\end{tabular}
\end{center}
\end{lemma}

The proof of Lemma \ref{lemma:setpartitionstat} is straightforward and involves expanding out the definition of $f(\pi,t)$, telescoping sums, and standard numerics. By Lemma \ref{lemma:setpartitionstat}, a set partition of $[k]$ into $s$ parts which minimizes $f(\pi,t)$ consists of one part of size $k-s+1$ and $s - 1$ singleton parts and hence 
\begin{equation}\label{eq:minpi}
 \min\{f(\pi,t) : |\pi| = s\} = (k-s+1) \binom{s-1}{t-1} + \binom{s-1}{t}.   
\end{equation}

In what follows, we denote a set partition of $[k]$ into $s$ parts which minimizes $f(\pi, t)$ by $\pi_s$. 

For $n$ large enough, and all $s$ where $t \leq s \leq k$, we will show \[ \frac{\sum_{ |\pi|=t} (1-t^{-t})^{Nf(\pi,t)}}{\sum_{ |\pi|=s} (1-t^{-t})^{Nf(\pi,t)}} \geq n^{s-t}. \]

Replacing the numerator with its largest term and each term in denominator with its largest term,  \[ \frac{\sum_{ |\pi|=t} (1-t^{-t})^{Nf(\pi,t)}}{\sum_{ |\pi|=s} (1-t^{-t})^{Nf(\pi,t)}}\geq \frac{(1-t^{-t})^{Nf(\pi_t,t)}}{S(k,s)(1-t^{-t})^{Nf(\pi_s,t)} } = \frac{1}{S(k,s)}(1-t^{-t})^{N(f(\pi_s,t)-f(\pi_t,t))}  \]
where $S(k,s)$ is the Stirling number of the second kind. Taking $n \geq S(k,s)$, a calculation using telescoping sums, Equation \eqref{eq:minpi}, and standard binomial coefficient bounds yields
\[ \frac{1}{S(k,s)}(1-t^{-t})^{N(f(\pi_s,t)-f(\pi_t,t))} \geq n^{s-t}. \] 
\begin{comment}
\begin{equation}\label{eq:equivcond}
k^{t-1} \frac{\alpha(s)-\alpha(t)}{\binom{k}{t}} \geq s-t+1 \Longrightarrow \frac{\rho(n,t)}{\rho(n,s)} \geq n^{s-t}.
\end{equation}

Observe that $\alpha(t+1)- \alpha(t) = (t-1)(k-t)$ and for $s>t$,
\begin{align*}
\alpha(s+1)-\alpha(s) &= \left[(k-s) \binom{s}{t-1} + \binom{s}{t}  \right] - \left[(k-s+1) \binom{s-1}{t-1} + \binom{s-1}{t}\right] \\
&=  (k-s) \left[ \binom{s}{t-1} - \binom{s-1}{t-1}  \right] + \left[\binom{s}{t}- \binom{s-1}{t-1} - \binom{s}{t-1}\right] = (k-s) \binom{s-1}{t-2}.
\end{align*}

\noindent Examining the left hand side in Equation \eqref{eq:equivcond} and using telescoping sums, \[ \frac{k^{t-1}}{\binom{k}{t}} (\alpha(s)-\alpha(t)) = \frac{k^{t-1}}{\binom{k}{t}} \sum_{i=t}^{s-1} (\alpha(i+1)-\alpha(i)) \geq \frac{k^{t-1}}{\binom{k}{t}} (s-t)(t-1)(k-t) \geq (s-t+1) \] since  $k^{t-1}(t-1)(k-t)\geq 2\binom{k}{t}$.
\end{comment}

Therefore, the index $s=t$ maximizes the right hand side of Equation \eqref{eq:expectationkt}, and hence \[ \mathbb E[U]  \leq (k-t+1) (n^t) \sum_{ |\pi|=t} (1-t^{-t})^{Nf(\pi,t)} < (k-t+1)n^t S(k,t) (1-t^{-t})^{N(k-t+1)} < 1 \] for our choice of $N$ provided $n\geq kS(k,t)$. Thus, \[ f_{k,t}(n) \leq \binom{k}{t}  \frac{(t+1)t^t \log_2 n}{(k-t+1) \log_2 e} =  \frac{(t+1)t^{t-1}}{\log_2 e} \binom{k}{t-1} \log_2n. \]

\subsection{Lower bound on $f_{k,2}(n)$}\label{subsec:lbfk2} 

In this section, we show
\begin{equation}\label{eq:nicebound}
 f_{k,2}(n) \geq \min\{m : {m \choose \lceil m/k \rceil} \geq n\}.
 \end{equation}

Let $\{H_1,H_2,\dots,H_m\}$ be a covering of $H_{k,2}(n)$ with $m = f_{k,2}(n)$ complete $k$-partite $k$-graphs. We recall
$H_{k,2}(n) = K_{n,n,\dots,n} \backslash M$, where $M$ is a perfect matching of $K_{n,n,\dots,n}$. For $i \in [k]$ and $j \in [n]$, define $A_{i,j} = \{H_r : x_{ij} \in V(H_r)\}$ and $\mathcal{A}_i = \{A_{i,j} : 1 \leq j \leq n\}$. As in (\ref{connection}), $(\mathcal{A}_1,\mathcal{A}_2,\dots,\mathcal{A}_k)$ is a Bollob\'{a}s $(k,2)$-tuple of size $n$. For convenience, for each $i \in [k]$, let $\phi_i:[k] \to [2]$ be so that $\phi_i^{-1}(1) = \{i\}$. 
Taking the sum of inequality from Theorem \ref{thm:genthreshk} with $t = 2$ over all $i \in [k]$,
\begin{equation}\label{eq:bsp_heccn}
\sum_{i = 1}^k \sum_{j=1}^{n} {|A_{1,j}(\phi_i) \cup A_{2,j}(\phi_i)| \choose |A_{1,j}(\phi_i)|}^{-1}  \leq k.
\end{equation}

We use this inequality to give a lower bound on $f_{k,2}(n) = m$. First we observe
\begin{equation}\label{eq:firstidentity}
 \sum_{r = 1}^m |V(H_r)| =  \sum_{j=1}^{n} \sum_{i = 1}^k |A_{i,j}|    =\sum_{j=1}^{n} \sum_{i = 1}^k |A_{1,j}(\phi_i)|.
\end{equation}

Let $\partial H$ denote the set of $(k - 1)$-tuples of vertices contained in some edge of a hypergraph $H$. Then
\begin{equation}\label{eq:secondidentity}
\sum_{r = 1}^m |\partial H_r \cap \partial M| = \sum_{j=1}^{n} \sum_{i = 1}^k |A_{2,j}(\phi_i)|.
\end{equation}

Putting the above identities together,
\begin{equation}\label{eq:thirdidentity}
\sum_{r = 1}^m |V(H_r)| + \sum_{r = 1}^m |\partial H_r \cap \partial M| = \sum_{j=1}^{n} \sum_{i = 1}^k (|A_{1,j}(\phi_i)| + |A_{2,j}(\phi_i)|).
\end{equation}

We note $|\partial H_r \cap \partial M| \leq |V(H_r)|/(k - 1)$, and therefore
\begin{equation}\label{inequality}
\sum_{r = 1}^m |\partial H_r \cap \partial M| \leq \frac{1}{k - 1} \sum_{r = 1}^m |V(H_r)|.
\end{equation}
It follows that
\begin{equation}\label{eq:inequality}
  \sum_{j=1}^{n} \sum_{i = 1}^k (|A_{1,j}(\phi_i)| + |A_{2,j}(\phi_i)|) \leq \frac{k}{k - 1} \sum_{r = 1}^m |V(H_r)|.
\end{equation}
Subject to the linear inequalities (\ref{eq:firstidentity}) and (\ref{eq:inequality}), the left side of (\ref{eq:bsp_heccn}) is minimized when $kn|A_{1,j}(\phi_i)| = \sum_{r = 1}^m |V(H_r)|$ and
 $kn(|A_{1,j}(\phi_i)| + |A_{2,j}(\phi_i)|) = (k - 1)|A_{1,j}(\phi_i)|$. Since $|V(H_r)| \leq (k - 1)n$ for all $r \in [m]$,
(\ref{eq:bsp_heccn}) implies ${m \choose \lceil m/k \rceil} \geq n$, which gives (\ref{eq:nicebound}). \qed

\subsection{Lower bound on $f_{k,k}(n)$}\label{subsection:lbfkk}
Let $\h = \{H_1,H_2,\dots,H_m\}$ be a minimal covering of $H_{k,k}(n)$ with complete $k$-partite $k$-graphs, so $m = f(H_{k,k}(n))$. Given a $k$-partite $k$-graph $H$, consider its $2$-shadow $\delta_2(H) = \{R \subset V(H) : |R|=k-2, R \subset e \text{ for some } e \in H  \}$. Let $\delta_2(\h) = \bigcup_{i=1}^m \delta_2(H_i)$.

Given $R \in \delta_2(\h)$ and $H_i \in \h$, let $H_i(R):=\{ e \in \binom{V(H_i)}{2}: e \cup R \in H_i \}$ be the possibly empty link graph of the edge $R$ in the hypergraph $H_i$ and let $V(H_i(R))$ be the set of vertices in the link graph. Observe that double counting yields 
\begin{equation}\label{eq:kminus2link}
\sum_{R \in \delta_2(\h)} \bigg( \sum_{i=1}^m |V(H_i(R))| \bigg) = \sum_{i=1}^m \bigg( \sum_{R \in \delta_2(H_i)} |V(H_i(R))| \bigg). 
\end{equation}

An optimization argument yields $|\delta_2(H_i)|$ is maximized when the parts of $H_i$ are of equal or nearly equal maximal size. Since $|V(H_i(R))|\leq 2(n-k+2)$, the right hand side of Equation \eqref{eq:kminus2link} is bounded above by 
\begin{equation}\label{eq:kminus2linkub}
\sum_{i=1}^m \bigg( \sum_{R \in \delta_2(H_i)} |V(H_i(R))| \bigg) \leq m \cdot \binom{k}{2} \cdot \left(\frac{n}{k}\right)^{k-2} \cdot 2(n-k+2).    
\end{equation}

For a lower bound on the left hand side of Equation \eqref{eq:kminus2link}, fix $R \in \delta_2(\h)$ and without loss of generality suppose that $R=\{ x_{1,1}, \ldots, x_{k-2,k-2} \}$. Let $Y= [k-1,n]$. Let $K_{Y,Y}$ be the complete bipartite graph with two distinct copies of $Y$ and $\M = \{ (x_{k-1,i}, x_{k,i} : i \in Y \}$ be a perfect matching in $K_{Y,Y}$. Then, $\{H_1(R), \ldots, H_m(R)\}$ forms a biclique cover of $K_{Y,Y} \setminus \M$. Applying the convexity result of Tarjan \cite[Lemma 5]{Tarjan}, \[\sum_{i=1}^m |V(H_i(R))| \geq (n-k+2) \log_2(n-k+2). \] 

Noting that $|\delta_2(\h)| =  \binom{k}{2}(n)_{(k-2)}$, the left hand side of Equation \eqref{eq:kminus2link} is bounded below by 
\begin{equation}\label{eq:kminus2linklb}
\sum_{R \in \delta_2(\h)} \bigg( \sum_{i=1}^m |V(H_i(R))| \bigg)  \geq \binom{k}{2}(n)_{(k-2)} (n-k+2) \log_2(n-k+2).   
\end{equation}

Comparing the bounds from Equation \eqref{eq:kminus2linkub} and Equation \eqref{eq:kminus2linklb}, \[ m  \geq \frac {(n)_{(k-2)}\log_2(n-k+2)}{  2\left(\frac{n}{k}\right)^{k-2} } \geq \frac{k^{k-2}}{2} \log_2 n \]
provided that $n$ is large enough. 

For $t \geq 3$ and $t<k$, the lower bound on $f_{k,t}(n)$ in Theorem \ref{thm:generalktbounds} is obtained from the lower bounds on $f_{t-1,t-1}(n-1)$ as follows: Let $\h = \{H_1,H_2,\dots,H_m\}$ be a minimal covering of $H_{k,t}(n)$ with complete $k$-partite $k$-graphs, so $m = f(H_{k,t}(n))$. Given $T \in \binom{[k]}{k-t+1}$, define $H_T \subset H_{k,t}(n)$ by \[ H_T:= \{  \{x_{1,i_1}, \ldots, x_{k,i_k} \} \in H_{k,t}(n) : i_j =1 \; \forall \; j \in T \}.   \] It follows that at least $f_{t-1,t-1}(n-1)$ of the complete $k$-partite $k$-graphs in $\h$ are needed to cover $H_T$. Moreover, for distinct $T,T' \in \binom{[k]}{k-t+1}$, the corresponding complete $k$-partite $k$-graphs from $\h$ are necessarily pairwise disjoint and hence \[ f_{k,t}(n) \geq \binom{k}{k-t+1} f_{t-1,t-1}(n-1) \geq \binom{k}{t-1} \frac{(t-1)^{t-3}}{2} \log_2 n \] provided that $n$ is large enough.   

\section{Concluding remarks}

$\bullet$ Our main theorem, Theorem \ref{thm:genthreshk} is tight for $t = 2$ and $k \geq 2$, as shown in Section \ref{subsec:sharpness}. It would be interesting to generalize this example to $2 < t \leq k$ to determine whether Theorem \ref{thm:genthreshk} is tight in general. The first open case is $t = k = 3$.

$\bullet$ A particular case of the Bollob\'{a}s set pairs inequality occurs when every set in $\mathcal{A}$ has size $a$ and every set in $\mathcal{B}$ has size $b$, and one obtains the tight bound $|\mathcal{A}| \leq {a + b \choose b}$. The generalization to Bollob\'{a}s $(k,t)$-tuples for $k \geq 3$ is equally interesting but wide open, as are potential generalizations to vector spaces -- see Lov\'{a}sz~\cite{L1,L2}.

$\bullet$  Orlin~\cite{Orlin} proved that the clique cover number $cc(K_n \backslash M)$ of a complete graph $K_n$ minus a perfect matching $M$ is precisely $\min\{m : 2{m - 1 \choose \lfloor m/2 \rfloor} \geq n\}$. Theorem \ref{thm:generalktbounds} yields lower bounds on the clique cover number of the complement of a perfect matching $M$ in the complete $k$-uniform hypergraph $K_n^k$:

\begin{cor}
Let $K_n^k\setminus M$ be the complement of a perfect matching in $K_n^k$. Then
\[  cc(K_n^k\setminus M) \geq \frac{\log_2\frac{n}{k}}{H(\frac{1}{k})} \geq \frac{k \log_2\frac{n}{k}}{\log_2 (ke)}. \]
\end{cor}

$\bullet$ It would be interesting to prove an analog of Equation \eqref{eq:nicebound} for $t\geq 3$. That is,  \[ f_{k,t}(n) \geq \min\{m : \binom{m}{\alpha_1, \ldots, \alpha_t}  \geq n_{(t-1)} \} \]
for some optimal $\alpha_1, \ldots, \alpha_t$. The difficulty here lies in determining effective bounds on $|A_{i,\sigma}(\phi)|$.

\end{document}